\documentclass[12pt, a4paper]{article}

\usepackage[utf8]{inputenc}
\usepackage{amsmath, amssymb, amsthm}
\usepackage{geometry}
\usepackage{xcolor}
\usepackage{hyperref}
\usepackage{graphicx}
\usepackage{booktabs}
\usepackage{tcolorbox}
\usepackage{listings}

\definecolor{darkblue}{rgb}{0.0, 0.2, 0.6}
\definecolor{codegreen}{rgb}{0,0.6,0}
\definecolor{codegray}{rgb}{0.5,0.5,0.5}
\definecolor{codepurple}{rgb}{0.58,0,0.82}
\definecolor{backcolour}{rgb}{0.95,0.95,0.92}

\hypersetup{colorlinks=true, linkcolor=darkblue, citecolor=darkblue, urlcolor=darkblue}

\lstdefinestyle{mystyle}{
    backgroundcolor=\color{backcolour},   
    commentstyle=\color{codegreen},
    keywordstyle=\color{magenta},
    numberstyle=\tiny\color{codegray},
    stringstyle=\color{codepurple},
    basicstyle=\ttfamily\footnotesize,
    breakatwhitespace=false,         
    breaklines=true,                 
    captionpos=b,                    
    keepspaces=true,                 
    numbers=left,                    
    numbersep=5pt,                  
    showspaces=false,                
    showstringspaces=false,
    showtabs=false,                  
    tabsize=4
}
\title{\textbf{A Formal Refutation of the Hypergeometric Parametric Extension for Reciprocal Binomial Sums}}
\author{Johar M. Ashfaque}
\date{\today}

\begin{document}

\maketitle

\begin{abstract}
Recent work by Pain \cite{pain2026} proposed a systematic approach to evaluating binomial sums involving reciprocals of binomial coefficients via Beta integrals. In particular, a parametric extension (Proposition 6.1) was introduced and claimed to admit a closed-form representation in terms of a terminating ${}_{2}F_{1}$ hypergeometric function. Through a combination of internal logical consistency checks, integral derivation analysis, and exact symbolic computation, we definitively prove that this parametric identity is false. 
\end{abstract}

\section{Introduction and Definitions}
Binomial sums involving reciprocals of binomial coefficients frequently appear in combinatorial analysis. The text under review \cite{pain2026} defines a parametric sum:
\begin{equation}
S_{n}(b,c;x) = \sum_{k=0}^{n}(-1)^{k}\binom{n}{k}\frac{x^{k}}{\binom{b+k}{c}} \label{eq:def}
\end{equation}
In Proposition 6.1 of the original work, the author claims that this sum evaluates to a specific hypergeometric finite expansion:
\begin{equation}
S_{n}(b,c;x) = \frac{c}{n+c}\frac{1}{\binom{n+b}{b-c}}{}_{2}F_{1}(-n,c+1;n+c+1;x) \label{eq:prop}
\end{equation}
While classical evaluations in the paper (such as Frisch's identity) are correct, the parametric generalization contains critical algebraic flaws.

\section{Logical Contradiction at $x=1$}
The most immediate disproof of Proposition 6.1 is that it fails to recover the verified base case of the paper. In Theorem 4.1 of \cite{pain2026}, the author successfully proves Frisch's identity, which evaluates the sum when the parameter $x = 1$:
\begin{equation}
S_{n}(b,c;1) = \frac{c}{n+c}\frac{1}{\binom{n+b}{b-c}}
\end{equation}
If Proposition 6.1 is valid, evaluating Eq.~\eqref{eq:prop} at $x = 1$ must logically yield Frisch's identity. This requires that the hypergeometric term evaluates to exactly 1. We test this using Gauss's Hypergeometric Theorem:
\begin{equation}
{}_{2}F_{1}(a,b;c_{param};1) = \frac{\Gamma(c_{param})\Gamma(c_{param}-a-b)}{\Gamma(c_{param}-a)\Gamma(c_{param}-b)}
\end{equation}
Substituting the author's parameters ($a = -n$, $b = c+1$, and $c_{param} = n+c+1$):
\begin{align*}
{}_{2}F_{1}(-n, c+1; n+c+1; 1) &= \frac{\Gamma(n+c+1)\Gamma((n+c+1) - (-n) - (c+1))}{\Gamma((n+c+1) - (-n))\Gamma((n+c+1) - (c+1))} \\
&= \frac{\Gamma(n+c+1)\Gamma(2n)}{\Gamma(2n+c+1)\Gamma(n)}
\end{align*}

\begin{tcolorbox}[colback=red!5,colframe=red!75!black,title=Contradiction]
The ratio $\frac{\Gamma(n+c+1)\Gamma(2n)}{\Gamma(2n+c+1)\Gamma(n)}$ does not equal 1. For example, if $n=2$ and $c=1$, the ratio evaluates to $\frac{36}{120} = \frac{3}{10}$. Therefore, the proposed generalization mathematically contradicts the author's own verified base case.
\end{tcolorbox}

\section{Calculus Autopsy: The Flawed Integral Derivation}
The root cause of the error lies in the evaluation of the Beta integral. The author starts with the correct parametric integral representation derived in Theorem 5.1 \cite{pain2026}:
\begin{equation}
S_{n}(b,c;x) = \int_{0}^{1}t^{c}(1-t)^{b-c}\Big[(b+1)(1-x(1-t))^{n}-nx(1-t)(1-x(1-t))^{n-1}\Big]dt
\end{equation}
The integrand explicitly consists of \textbf{two} distinct terms inside the brackets. However, in the proof of Proposition 6.1:
\begin{enumerate}
    \item \textbf{Omission of the Derivative Term:} The author drops the second term of the integrand, $-nx(1-t)(1-x(1-t))^{n-1}$, without any algebraic justification, proceeding to integrate only the first term.
    \item \textbf{Fabricated Factorials:} Even when strictly integrating the first term by expanding $(1-x(1-t))^n$ via the binomial theorem, the resulting Beta integrals yield coefficients involving $\frac{c!(b-c+j)!}{(b+j+1)!}$. There is no valid algebraic pathway to transform these factorials into the specific Pochhammer ratio $\frac{(c+1)_j}{(n+c+1)_j}$ required to form the claimed ${}_{2}F_{1}$ series.
\end{enumerate}

\section{Exact Symbolic Verification}
To preclude any arguments regarding numerical floating-point errors, we provide a Python script utilizing the \texttt{sympy} library. This script treats $x$ as an abstract algebraic symbol, expanding both the definitional summation and the author's proposed formula into exact polynomials.

\begin{lstlisting}[language=Python, caption=Exact Symbolic Verification using SymPy]
import sympy as sp

def symbolic_proof(n_val, b_val, c_val):
    x = sp.Symbol('x')
    
    # 1. Calculate the EXACT LHS (Definition of the sum)
    lhs_sum = 0
    for k in range(n_val + 1):
        binom_nk = sp.binomial(n_val, k)
        binom_bk_c = sp.binomial(b_val + k, c_val)
        term = ((-1)**k) * binom_nk * (x**k) / binom_bk_c
        lhs_sum += term
        
    lhs_poly = sp.simplify(lhs_sum)
    
    # 2. Calculate the EXACT RHS (Proposition 6.1)
    prefactor = sp.Rational(c_val, n_val + c_val) * (1 / sp.binomial(n_val + b_val, b_val - c_val))
    
    rhs_series = 0
    for j in range(n_val + 1):
        binom_nj = sp.binomial(n_val, j)
        poch_num = sp.factorial(c_val + 1 + j - 1) / sp.factorial(c_val)
        poch_den = sp.factorial(n_val + c_val + 1 + j - 1) / sp.factorial(n_val + c_val)
        
        term = ((-1)**j) * binom_nj * sp.Rational(poch_num, poch_den) * (x**j)
        rhs_series += term
        
    rhs_poly = sp.simplify(prefactor * rhs_series)
    
    # 3. Check for algebraic equivalence
    is_equivalent = sp.simplify(lhs_poly - rhs_poly) == 0
    return sp.expand(lhs_poly), sp.expand(rhs_poly), is_equivalent

# Test with n=2, b=3, c=1
lhs, rhs, match = symbolic_proof(2, 3, 1)
print(f"LHS Polynomial: {lhs}")
print(f"RHS Polynomial: {rhs}")
print(f"Match: {match}")
\end{lstlisting}

\subsection{Execution Results}
When executed for parameters $n=2$, $b=3$, and $c=1$, the script outputs exact fractional polynomials:
\begin{itemize}
    \item \textbf{LHS (True Sum):} $\frac{1}{5}x^{2} - \frac{1}{2}x + \frac{1}{3}$
    \item \textbf{RHS (Claimed Formula):} $\frac{1}{100}x^{2} - \frac{1}{30}x + \frac{1}{30}$
\end{itemize}
Because these two polynomials are algebraically distinct, it is a definitive mathematical fact that Proposition 6.1 does not hold.

\section{Conclusion}
Through logical base-case evaluation, an autopsy of the integral manipulation, and exact symbolic computation, we have conclusively refuted the hypergeometric parametric extension proposed in Proposition 6.1 of \cite{pain2026}. The closed-form representation fails due to the improper handling of the Beta integral expansion.

\end{document}